\documentclass[11pt,oneside]{amsart}

\usepackage[usenames,dvipsnames,svgnames,table]{xcolor}
\usepackage[colorlinks=true, pdfstartview=FitV, linkcolor=blue, citecolor=blue, urlcolor=darkblue]{hyperref}

\usepackage{microtype}
\usepackage{geometry}                % See geometry.pdf to learn the layout options. There are lots.
\usepackage{graphicx}
\usepackage{mathrsfs}
\usepackage{amssymb}
\usepackage{amsfonts}
\usepackage{mathrsfs}
\usepackage{epstopdf}
\usepackage{lscape}
\usepackage[utf8]{inputenc}
\usepackage{tikz,caption}
\usepackage{enumitem}
\setlist[itemize]{leftmargin=2em}
\setlist[enumerate]{leftmargin=2em}
\usepackage{booktabs}
\usepackage{multirow}
\usepackage{mathtools}
\usepackage[linesnumbered,ruled]{algorithm2e}

\definecolor{darkblue}{rgb}{0.0,0,0.7} % darkblue color
\definecolor{darkred}{rgb}{0.7,0,0} % darkred color
\definecolor{darkgreen}{rgb}{0, .6, 0} % darkgreen color

\newcommand{\defncolor}{\color{darkred}}
\newcommand{\defn}[1]{{\defncolor\emph{#1}}} % emphasis of a definition

\newtheorem{theorem}{Theorem}[section]
\newtheorem{prop}[theorem]{Proposition}
\newtheorem{cor}[theorem]{Corollary}
\newtheorem{lemma}[theorem]{Lemma}

\theoremstyle{definition}

\newtheorem{example}[theorem]{Example}
\newtheorem{remark}[theorem]{Remark}
\numberwithin{equation}{section}

\usepackage[colorinlistoftodos]{todonotes}

\newcommand{\qbinom}[2]{\left[ \begin{smallmatrix}#1\\#2\end{smallmatrix} \right]}

\title{Quasisymmetric harmonics of the exterior algebra}
\author{Nantel Bergeron,
Kelvin Chan,
Farhad Soltani,
Mike Zabrocki}
\address{Dept. of Math. and Stat.\\ York  University\\ To\-ron\-to, Ontario M3J 1P3\\ CANADA}
\email{bergeron@yorku.ca}
\email{ktychan@yorku.ca}
\email{farhad.soltani91@gmail.com}
\email{zabrocki@yorku.ca}

\thanks{This work is supported in part by York Research Chair and NSERC.
This paper originated in a working session at the Algebraic
Combinatorics Seminar at York University}
\date{May 28, 2022}



\begin{document}

\begin{abstract}
We study the ring of quasisymmetric polynomials in $n$ anticommuting  (fermionic) variables.
Let $R_n$ denote the polynomials in $n$ anticommuting variables. The main results of this paper show the following interesting facts about  quasisymmetric polynomials in anticommuting variables:
\begin{enumerate}
\item The quasisymmetric polynomials in $R_n$ form a commutative sub-algebra of $R_n$.
\item There is a basis of the quotient of $R_n$ by the ideal $I_n$ generated by the
quasisymmetric polynomials in $R_n$ that is indexed by ballot sequences.
The Hilbert series of the quotient is given by
 $$ \text{Hilb}_{R_n/I_n}(q) = \sum_{k=0}^{\lfloor{n/2}\rfloor} f^{(n-k,k)} q^k\,,$$
 where $f^{(n-k,k)}$ is the  number  of standard tableaux of shape $(n-k,k)$.
\item There is a basis of the ideal generated by quasisymmetric polynomials
that is indexed by sequences that break the ballot condition
\end{enumerate}
\end{abstract} 

\maketitle

%%%%%%%%%%%%%%%%%%%%%%%%%%%%%%%%%%%%%%
%%%%%%%%%%%%%%%%%%%%%%%%%%%%%%%%%%%%%%
%%%%%%%%%%%%%%%%%%%%%%%%%%%%%%%%%%%%%%
\section{Introduction}
The study of coinvariants of groups dates back to Shephard-Todd and Chevalley~\cite{ST,Chevalley} and has fruitfully produced many connections between algebra, combinatorics and physics.  Motivated by recent developments in coinvariants of symmetric groups and symmetric functions theory incorporating fermionic variables, we study a coinvariant-like quotient of an exterior algebra obtained by the quotient of the ideal generated by quasisymmetric functions in fermionic variables.  The quotient has a dimension that can be interpreted as the number of ballot sequences (or other interpretations, see for instance the OEIS \cite{OEIS} sequences \href{https://oeis.org/A008315}{A008315} and \href{https://oeis.org/A001405}{A001405}).

A notable feature of many quotients similar to coinvariants is their amenability to combinatorial methods.  One well-known example is the coinvariant ring of the symmetric group.  It is the quotient of the polynomial ring $\mathbb{Q}[x_{1},\dots,x_{n}]$ in commuting variables by the ideal generated by the symmetric polynomials with no constant term.  As an $\mathcal{S}_{n}$ representation, this quotient is naturally graded and is well-known to be isomorphic to the regular representation.  Many useful bases of this space have been found by studying combinatorics related to permutations.  For more details, see the nice surveys of \cite{B, GH, MacSchub, Manivel}.

This line of inquiry inspired Garsia and Haiman \cite{GH96,H} to consider the ring of diagonal harmonics, a similar quotient in two sets of commuting variables as an $\mathcal{S}_{n}$ module.  Haiman's work \cite{H2} showed that the diagonal harmonics have a deep connection to the theory of Macdonald polynomials.  A combinatorial expression for the Frobenius image of the diagonal harmonics known as the Shuffle Conjecture \cite{HHLRU} showed that the module structure is closely related to the combinatorics of parking functions and can be described in terms of certain labelled Catalan paths.  This connection relating the symmetric functions and the combinatorial expression was proven in \cite{CM} and is now known as the Shuffle Theorem.

The connection between the combinatorics and the symmetric function expressions of the Shuffle Theorem have been generalized \cite{HRW} and proven \cite{DM} to an expression known as the Delta Conjecture.  The last author with the group at the Fields Institute \cite{Z} proposed a deformation of diagonal harmonics to two sets of commuting variables and one set of anticommuting variables. In this case, the connection of representation theoretic interpretation to the symmetric function expression remains open.  The symmetric function expressions and representation theoretic interpretation was extended further to include the quotient of two sets of commuting and two sets of anticommuting variables in \cite{DIW} to what is known as the Theta Conjecture.  At present, this also remains an open conjecture, but progress has been made on some special cases \cite{IRR, KR, SW, SW2}.

The ring of quasisymmetric polynomials $QSym_{n}$ contains the ring of symmetric polynomials $Sym_{n}$.  Many combinatorial structures of $QSym_{n}$ parallel that of $Sym_{n}$.  Hivert described a Temperley-Lieb $TL_{n}$ action on $\mathbb{Q}[x_{1},\dots,x_{n}]$ making $QSym_{n}$ exactly its trivial representation \cite{Hi}.  In 2003, Aval, F. Bergeron, and the first author studied $QSym$ coinvariant spaces obtained by replacing the ideal of non-constant symmetric functions with the ideal of non-constant quasisymmetric functions \cite{AB,ABB}.  Surprisingly they found that dimensions of $QSym$ coinvariants are equal to the Catalan numbers.  At the heart of their argument is a recursion built from Catalan paths.  Li extended this argument to study some components of coinvariant spaces of diagonally quasisymmetric functions \cite{L}.

Motivated by physics, Desrosiers, Lapointe, and Mathieu \cite{DLM,DLM2} introduced symmetric functions with one set of commuting and one set of anticommuting variables known as symmetric function in superspace.
The commuting variables encode bosons while the anticommuting ones encode fermions, hence the anticommuting variables are sometimes referred to as ``fermionic variables.''  The Hopf algebra structure of the ring of symmetric functions in superspace was extended to quasisymmetric functions in superspace \cite{FLP} and so a natural question is to extend the study of coinvariants of polynomial rings with commuting and anticommuting variables to the quotients of  these polynomial rings by the ideal generated by ``super'' quasisymmetric polynomials.

Parallel to the Delta Conjecture or Theta Conjecture, one ideally would like to understand quasisymmetric coinvariants in multiple sets of commuting and anticommuting variables.  Our study of quasisymmetric coinvariant spaces in one set of anticommuting variables is a first step in that study.  We denote polynomials in anticommuting variables by $R_n$. The main results of this paper show the following interesting facts about symmetric and quasisymmetric functions in anticommuting variables:
\begin{enumerate}
\item The quasisymmetric polynomials in $R_n$ form a commutative sub-algebra of $R_n$
(Proposition \ref{prop:comm}).
\item That $R_n$ is free over the ring of symmetric polynomials (Proposition \ref{prop:free}).
\item There is a basis of the quotient of $R_n$ by the ideal $I_n$ generated by the
quasisymmetric polynomials in $R_n$ that is indexed by ballot sequences (Proposition \ref{prop:harmbasis}).
The Hilbert series of the quotient is given by
 $$ \text{Hilb}_{R_n/I_n}(q) = \sum_{k=0}^{\lfloor{n/2}\rfloor} f^{(n-k,k)} q^k\,,$$
 where $f^{(n-k,k)}$ is the  number  of standard tableaux of shape $(n-k,k)$ (Corollary~\ref{cor:hilb}).
\item There is a basis of the ideal generated by quasisymmetric polynomials
that is indexed by sequences that break the ballot condition
(Theorem \ref{thm:basisofideal}) and a minimal Gr\"obner basis
that is a subset of this basis (Corollary \ref{cor:minimalGB}).
\end{enumerate}

\subsection{Acknowledgement} We are grateful to Yohana Solomon for her participation and discussions during this project. We also thank D. Grinberg for a careful reading of the menuscript and helpful comments.

%%%%%%%%%%%%%%%%%%%%%%%%%%%%%%%%%%%%%%
%%%%%%%%%%%%%%%%%%%%%%%%%%%%%%%%%%%%%%
%%%%%%%%%%%%%%%%%%%%%%%%%%%%%%%%%%%%%%
\section{Quasisymmetric invariants on the exterior algebra}

Fix $n$ a positive integer and
let $R_n = {\mathbb Q}[\theta_1, \theta_2, \ldots, \theta_n]$ be the
polynomial ring in anticommuting variables.
The ring $R_n$ is isomorphic to the exterior algebra of a vector
space of dimension $n$.  The variables of this ring satisfy the relations
\[
\theta_i \theta_j = - \theta_j \theta_i \hbox{ if } 1 \leq i \neq j \leq n
\qquad\hbox{and}\qquad \theta_i^2 = 0 \hbox{ for }1 \leq i \leq n~.
\]
Since these conditions impose that a monomial in $R_n$ has no repeated variables,
the monomials are in bijection with subsets of $\{1,2,\ldots, n\}$
and the dimension of $R_n$ is therefore equal to $2^n$.

Denote $[n] := \{1,2, \ldots,n\}$ and
let $A = \{a_1 < a_2 < \cdots < a_r \} \subseteq [n]$.
We define $\theta_A := \theta_{a_1} \theta_{a_2} \cdots \theta_{a_r}$,
then the set of monomials $\{ \theta_A \}_{A \subseteq [n]}$ is a basis of $R_n$.

We define an action on monomials of $R_n$ and extend this action linearly.
For each integer $1 \leq i < n$, let $\pi_i$ be an operator on $R_n$
that is defined by
\begin{equation}\label{eq:pi}
\pi_i(\theta_A) = \begin{cases}
\theta_{A} & \hbox{ if } i, i+1 \in A\hbox{ or }i, i+1 \notin A\\
\theta_{A \cup \{i+1\} \backslash \{i\}} & \hbox{ if } i\in A\hbox{ and }i+1 \notin A\\
\theta_{A \cup \{i\} \backslash \{i+1\}} & \hbox{ if } i+1\in A\hbox{ and }i \notin A
\end{cases}~.
\end{equation}
These operators instead of exchanging an $i$ for an $i+1$ like the symmetric group
action, have the effect of shifting the indices of the variables (if possible).  They
are therefore known as quasisymmetric operators.  They were studied in depth by
Hivert \cite{Hi}.  The operators are not multiplicative on $R_n$ in general since, for example,
\[
\pi_1( \theta_{1} \theta_{2})
= \theta_1 \theta_2
= - \pi_1( \theta_{1}) \pi_1(\theta_{2})~.
\]
They are also not multiplicative when they act on the polynomial ring
in commuting variables.

A polynomial that is invariant under the action of quasisymmetric operators
is said to be quasisymmetric invariant (or just `quasisymmetric').
The quasisymmetric invariants of $R_n$ are
linearly spanned by the elements:
\begin{equation}\label{eq:defF}
F_{1^r}(\theta_1, \theta_2, \ldots, \theta_n) := \sum_{\substack{A \subseteq [n]\\|A|=r}} \theta_A~.
\end{equation}
The symbols $F_{1^r}$ for the elements borrows the notation from the
polynomial ring in commuting variable invariants known as the `fundamental
quasisymmetric polynomials.'  The commuting polynomial quasisymmetric
invariants are indexed by compositions.

\begin{remark}
As expressing polynomials with listing the variables
(e.g. $p(\theta_1, \theta_2, \ldots, \theta_n)$) can be notational cumbersome,
there will be points where we will drop the variables in the expressions
and this will indicate that the polynomials are in the
variables $\theta_1, \theta_2, \ldots, \theta_n$.  There will also
be expressions where some polynomials have fewer variables and there
we will indicate this by listing the variables.
\end{remark}

%%%%%%%%%%%%%%%%%%%%%%%%%%%%%%%%%%%%%%
\subsection{Quasisymmetric functions generate a commutative subalgebra}
In \cite{FLP}, the authors showed that the quasisymmetric functions in
one set of commuting variables and one set of anticommuting variables
forms a graded Hopf algebra.  This implies that the quasisymmetric functions
in one set of anticommuting variables are closed under multiplication
and the space is spanned by one element at each non-negative degree.
It follows that for $r, s \geq0$ that there exists a (possibly $0$)
coefficient $a_{r,s}$ such that
\begin{equation}\label{eq:qsalg}
F_{1^r} F_{1^s} = a_{r,s} F_{1^{r+s}}\,.
\end{equation}
If $r+s>n$ then $F_{1^{r+s}} = 0$ by definition and so the only relevant coefficient $a_{r,s}$ is when
$r+s \leq n$.

\begin{remark}
In  the notation of \cite{FLP}, $F_{1^r}=M_{\dot{0}^r}=L_{\dot{0}^r}$ where $\dot{0}^r=(\dot{0},\dot{0},\ldots,\dot{0})$ a composition of length $r$.
The fermionic degree of $F_{1^r}$ is exactly $r$.
 In~\cite{FLP}, they show that $a_{r,s}$ exists and express it as a sum of $\pm 1$, but they do not give an explicit formula.
Furthermore, they indicate  that $a_{r,s}=(-1)^{rs}a_{s,r}$. Here we shall compute exactly $a_{r,s}$
and the formula shows the subalgebra generated by the $F_{1^r}$ is commutative.
\end{remark}

\begin{prop}\label{prop:comm}
The constants $a_{r,s}$ in Equation~\eqref{eq:qsalg} satisfy the following equation.
$$ a_{r,s}=
\begin{cases}
	0  &\text{if $r,s$ are both odd,}\\
	\\
	\left({\lfloor \frac{r+s}{2} \rfloor \atop \lfloor \frac{r}{2} \rfloor}\right)&\text{otherwise.}
\end{cases}
$$
\end{prop}
A remark brought to our attention by D. Grinberg \cite{DG} shows that $a_{r,s}$ is equal to the $q$-binomial coefficient $\qbinom{r+s}{r}_q$ evaluated at $q \rightarrow -1$ \cite[Equation~(185) on page 291]{DG2}.

\begin{proof}
For completeness, we give a proof not assuming any results of~\cite{FLP}. Using Equation~\eqref{eq:defF} we have
$$F_{1^r}F_{1^s}= \sum_{\substack{A \subseteq [n]\\|A|=r}}  \sum_{\substack{B \subseteq [n]\\|B|=s}} \theta_A \theta_B
= \sum_{\substack{C \subseteq [n]\\|C|=r+s}} \Big(\sum_{\substack{A\subseteq C \\|A|=r}} (-1)^{|\{b<a\,|\,a\in A,\, b\in C\setminus A\}|} \Big) \theta_C\,.
$$
To see the second equality, we remark that the product $\theta_A\theta_B=0$ if $A\cap B\ne \emptyset$.
Furthermore, if $A\cap B= \emptyset$, then for $C=A\cup B$
we have $B=C\setminus A$ and  $\theta_A\theta_B=(-1)^{|\{b<a\,|\,a\in A, b\in C\setminus A\}|} \theta_C$,
where the sign is the number of interchanges needed to sort $A$ followed by $B$ into $C$.
This does not depend on the values of the elements of $C$, but only on how $A$ is chosen inside $C$.
This shows that we get the same coefficient for all $C$ of size $r+s$ and therefore
$F_{1^r} F_{1^s} = a_{r,s} F_{1^{r+s}}$ with
\begin{equation}\label{eq:signa}
a_{r,s}=\sum_{\substack{A\subseteq \{1,2,\ldots, r+s\}\\ |A|=r}}
(-1)^{|\{1\le b<a\le r+s\,|\,a\in A, b\not\in A\}|}
\end{equation}
by choosing $C=\{1,2,\ldots,r+s\}$.

Let ${C\choose r}=\{A\subseteq C, |A|=r\}$. We define a sign-reversing involution $\Phi\colon {C\choose r}\to{C\choose r}$ as follows.
For $A\in  {C\choose r}$, let $\gamma(A)=\gamma_1\gamma_2\cdots\gamma_{r+s}\in\{0,1\}^{r+s}$, be the sequence such that
 $\gamma_i=1$ if $i\in A$, and $\gamma_i=0$ otherwise.
 We look at the  entries of  $\gamma(A)$ two by two and find the smallest $j$ (if it exists) such that the  pair
 $\gamma_{2j-1}\gamma_{2j}$ is not $00$ or $11$. If there is no such pair, we let $\Phi(A)=A$.   If we find such pair we define the involution $\Phi(A)=A'$ where $A'$ is such that $\gamma(A')$ is obtained from $\gamma(A)$ by interchanging $01\leftrightarrow 10$ in position $2j-1,2j$. If $r$ and $s$ are both odd, then there must be at least one occurrence of $01$ or $10$ and there are no fixed points of this involution.

We let
 $$Inv(A)=\{1\le b<a\le r+s\,|\,a\in A, b\not\in A\}=\{1\le \ell<t\le r+s\,|\,\gamma_\ell=0,\gamma_t=1 \}\,$$
 where $\gamma(A)=\gamma_1\gamma_2\cdots\gamma_{r+s}$.
As long as $(t,\ell)\ne(2j-1,2j)$ there is a bijection between $(t,\ell)\in Inv(A)$ and
$(t',\ell)\in Inv(A')$ interchanging the 1 and 0 in positions $2j-1$ and $2j$.
The pair $(2j-1,2j)$ is in only one of $Inv(A)$ or $Inv(A')$ but not the other. Therefore
$$(-1)^{|\{1\le b<a\le r+s\,|\,a\in A, b\not\in A\}|} = -(-1)^{|\{1\le b<a\le r+s\,|\,a\in A', b\not\in A'\}|}\,.$$
If $\Phi(A)=A$, we have that $|Inv(A)|$ is even since we can match the pairs two-by-two.
If $r$ is odd and $s$ is even, then the only $A \in \binom{C}{r}$ have
$r+s \in A$ and $|Inv(A)| = |Inv(A \backslash \{ r+s \})| + s$.
Therefore $\Phi$ is a sign reversing involution and all fixed points contribute
in Equation~\eqref{eq:signa} with a $+1$. Therefore
$$a_{r,s}=\Big|\big\{ A\in {C\choose r}\big| \Phi(A)=A\big\}\Big|=\left({\lfloor \frac{r+s}{2} \rfloor \atop \lfloor \frac{r}{2} \rfloor}\right)\,,$$
since there are a total of $\lfloor \frac{r+s}{2} \rfloor$ pairs $2j-1,2j$
in a sequence of length $r+s$ and we must have $\lfloor \frac{r}{2} \rfloor$ of them equal to $11$
and all others equal to $00$ in order to get $\Phi(A)=A$.
\end{proof}

The generating series for the coefficients $F(x,y) = \sum_{r,s \ge 0} a_{r,s} x^{r}y^{s}$ is
equal to $\frac{1 + x + y}{1 - x^{2} - y^{2}}$
and the OEIS \cite{OEIS} sequence number is \href{https://oeis.org/A051159}{A051159}.
This can be derived from Proposition \ref{prop:comm} using standard techniques of generating functions. 
%\mike{TODO: add comment to OEIS sequence}

One consequence of Proposition~\ref{prop:comm} is that $a_{r,s}=a_{s,r}$ for all $r,s\ge 0$. Remark that this does not contradict the fermionic law stating that $a_{r,s}=(-1)^{rs}a_{s,r}$ since $a_{r,s}=0$ when both $r,s$ are odd. Therefore we have shown the following corollary.

\begin{cor}
The subalgebra generated by quasisymmetric invariants $\{F_{1^r}|r\ge 0\}$ is commutative.\hfill$\square$
\end{cor}

\subsection{The ideal generated by symmetric invariants}
The symmetric invariants ${\rm Sym}_{R_n}$ of $R_n$ are very small since a basis consists of only two
elements $1$ and $F_1(\theta_1, \theta_2, \ldots, \theta_n)$.  Therefore the ideal generated by the invariants of non-zero
degree, which we shall denote $J_n$, is generated by a single element
$F_1$.
We begin by considering the symmetric
coinvariants of $R_n$, the quotient ring $R_n/J_n$.
Because the ideal $J_n$ is principal we can
understand this quotient with much more detail.
This quotient ring is a special case of the ring
recently studied in \cite{IRR,KR}.

Recall that ${\rm dim}~R_n = 2^n$,
and if we consider the quotient $R_n/J_n$ it is isomorphic to $R_{n-1}$ since
in this algebra $\theta_n = - \theta_1 -\theta_2 - \cdots - \theta_{n-1}$.
Let $A \subseteq [n-1]$ and $A' = A \cup \{n\}$, then
the map which sends
$\theta_{A'}$ to
$$-\theta_{A}(\theta_1 + \theta_2 + \cdots + \theta_{n-1}) \otimes 1
+ \theta_{A} \otimes F_1\quad \in\ R_n/J_n \otimes {\rm Sym}_{R_n}$$
and $\theta_{A}$ to
$$\theta_{A} \otimes 1 \quad \in \ R_n/J_n \otimes {\rm Sym}_{R_n}$$
is an algebra isomorphism.  Since this map describes the image
for each monomial in $R_n$, we have the following proposition.

\begin{prop} \label{prop:free}
For each $n \geq 1$,
$$R_n \cong R_n/J_n \otimes {\rm Sym}_{R_n},$$
as an algebra.  That is, $R_n$ is free over ${\rm Sym}_{R_n}$.
\end{prop}

%%%%%%%%%%%%%%%%%%%%%%%%%%%%%%%%%%%%%%
\subsection{The ideal generated by the quasisymmetric invariants}

Define an ideal of $R_n$ generated by the quasisymmetric invariants as
\[
I_n := \left< F_{1^r}(\theta_1, \theta_2, \ldots, \theta_n) : 1 \leq r \leq n \right>.
\]

\begin{remark}
Note that since the operators $\pi_i$ are not multiplicative, it
is unlikely to be the case that $I_n$ as an ideal is invariant
under the action of the $\pi_i$.  Indeed, we find that for $n=4$,
\[
\theta_2 F_{1}(\theta_1, \theta_2, \theta_3, \theta_4) =
-\theta_1 \theta_2 + \theta_2 \theta_3 + \theta_2 \theta_4.
\]
If we apply $\pi_1$ to this element, we obtain
\[
\pi_1(\theta_2 F_{1}(\theta_1, \theta_2, \theta_3, \theta_4)) =
-\theta_1 \theta_2 + \theta_1 \theta_3 + \theta_1 \theta_4.
\]
This is not in $I_4$.
\end{remark}

The \emph{exterior quasisymmetric coinvariants}\footnote{We borrow
the name `coinvariant' space even though the generators, and not the whole ideal, is
invariant under the quasisymmetric operators.} are defined to be
\[
EQC_n := R_n/I_n~.
\]

%%%%%%%%%%%%%%%%%%%%%%%%%%%%%%%%%%%%%%
\subsection{Differential operators on the exerior algebra}\label{ssec:harm}
We can define a set of differential operators on $R_n$ which
will permit us to define the orthogonal complement to the
ideal and a notion of quasisymmetric harmonics.

The operators $\partial_{\theta_i}$ act on monomials in $R_n$
by
\[
\partial_{\theta_i}( \theta_A ) = \begin{cases}
(-1)^{\#\{ j \in A: j<i\}}\theta_{A \backslash \{i\}}&\hbox{ if }i \in A\\
0&\hbox{ if }i \notin A
\end{cases}~.
\]

The operators can equally be characterized by the action that $\partial_{\theta_i}(1) = 0$
and the commutation relations
\[
\partial_{\theta_i} \partial_{\theta_j}=-\partial_{\theta_j} \partial_{\theta_i}
\hbox{ if } 1 \leq i \neq j \leq n
\qquad\hbox{and}\qquad
\partial_{\theta_i}^2 = 0\hbox{ for }1 \leq i \leq n
\]
\[
\partial_{\theta_i} \theta_j=-\theta_j \partial_{\theta_i}
\hbox{ if } 1 \leq i \neq j \leq n
\qquad\hbox{and}\qquad
\partial_{\theta_i} \theta_i = 1\hbox{ for }1 \leq i \leq n~.
\]

For a monomial $\theta_A = \theta_{a_1} \theta_{a_2} \cdots \theta_{a_r}$,
let $\overline{\theta_A} = \theta_{a_r} \theta_{a_{r-1}} \cdots \theta_{a_1}$ represent
reversing the order of the variables in the monomial.  Extend this notation to both
differential operators and polynomials (and polynomials of differential operators)
by extending the notation linearly.

We can define an inner product on $R_n$ by setting for $p,q \in R_n$.
\[
\left< p, q \right> = \overline{p(\partial_{\theta_1}, \partial_{\theta_2}, \ldots, \partial_{\theta_n})}
q( \theta_1, \theta_2, \ldots, \theta_n)|_{\theta_1=\theta_2 = \cdots=\theta_n=0}~.
\]
The monomials of $R_n$ form an orthonormal basis of the space with respect to this
inner product.

Define the orthogonal complement to $I_n$ with respect to the inner product as
the set
\begin{align}
EQH_n :&= \left\{ q \in R_n : \left< p, q \right> = 0 \hbox{ for all } p \in I_n \right\}
\label{eq:set1}\\
  &=\left\{ q \in R_n : p(\partial_{\theta_1}, \partial_{\theta_2}, \ldots, \partial_{\theta_n})
q= 0 \hbox{ for all } p \in I_n \right\}~.\label{eq:diffeqs}
\end{align}
The second equality follows from the fact that $I_n$ is an ideal
and shows that $EQH_n$ is also the solution space of
a system of differential equations.  The containment of the set
in Equation \eqref{eq:diffeqs} inside the set in Equation \eqref{eq:set1} is clear.
For the reverse inclusion, take an element $q$ which is not in the
set in Equation \eqref{eq:diffeqs}, we assume
for some $p \in I_n$ that $p(\partial_\theta) q = c \theta^\alpha$ plus possibly some other terms,
but then $\overline{p \theta^\alpha} \in I_n$
and $\left< \overline{p \theta^\alpha}, q \right> = c$ which implies that $q$ is not in the
the set in Equation \eqref{eq:set1}.
%\mike{Let $q$ be such that $\left<p,q\right>=0$ for all $p \in I_n$.
%Assume (by contradiction) that $\overline{p(\partial)}q \neq 0$ and
%is instead $c x^\alpha$ + terms of lower degree with $c\neq0$. But then
%$\left< p \overline{x^\alpha}, q \right>=c$ which is a contradiction.}
We refer to $EQH_n$ as the \emph{exterior quasisymmetric harmonics}.\footnote{
The harmonics and diagonal harmonics borrows the name from the physics literature
because the harmonic operator $\partial_1^2 + \partial_2^2 + \cdots + \partial_n^2$
is symmetric in the differential operators.  In the case of the exterior algebra,
this operator acts as zero and yet we persist by borrowing the name from the
analogous spaces of commuting variables.
}

It is clear that the monomials of $R_n$ form an orthonormal basis of the space
with respect to the inner product, hence the inner product is positive definite.
It follows that since $EQH_n$ is the orthogonal complement of the ideal $I_n$ in $R_n$,
then the following result must hold.
\begin{prop} \label{prop:EQC_EQH} For all $n \geq 1$, as graded vector spaces
\[
EQC_n \simeq EQH_n~.
\]
\end{prop}

We will conclude this section  by constructing a set of linearly independent elements inside $EQH_n$, which will
give us a lower  bound on the dimension of $EQC_n$. In Section~\ref{sec:ballotbasis} we will see that this is also an upper bound,
thus concluding that our  set is in fact a basis. To compute $EQH_n$ we need to solve the differential equations in Equation~\eqref{eq:diffeqs}.
Remark first  that since $I_n$ is an ideal, we do not need to take all $p\in I_n$, but it is enough to solve for the generators $p=F_{1^r}$ for  $1\le r\le n$.
We can reduce that further using Proposition~\ref{prop:comm} as noted in the following lemma.
\begin{lemma}\label{lem:idealgen}
For $n\ge 2$ we have $I_n$ is the ideal generated by $F_1$ and $F_{1^2}$.
\end{lemma}

\begin{proof} Clearly we have that the ideal generated by $F_1, F_{1^2}$ is contained in $I_n$.
For the converse we note that for each $k\geq 1$
there are non-zero coefficients $a$ and $a'$ such that
\[
a F_{1^{2k}} = (F_{1^2})^k \hskip .2in\hbox{   and   }\hskip .2in a' F_{1^{2k+1}} = (F_{1^2})^k F_1
\]
hence all of the generators of $I_n$ are contained in the ideal generated by $F_1, F_{1^2}$.
%For the  converse, consider first $2\le r=2k\le n$. Using Proposition~\ref{prop:comm}, we have that
%$aF_{1^r} = \big(F_{1^2}\big)^k$ for some $a\ne 0$. Therefor $F_{1^r} \in \langle F_1, F_{1^2}  \rangle$. For $2\le r=2k+1\le n$, we have $aF_{1^r} = \big(F_{1^2}\big)^kF_1$, which again show $F_{1^r} \in \langle F_1, F_{1^2}  \rangle$. We conclude that $I_n\subseteq  \langle F_1, F_{1^2}  \rangle$.
\end{proof}
From this we conclude that
\begin{equation}\label{eq:defEQH}
EQH_n =  \Big\{ q \in R_n :  \quad\sum_{1\le i\le n} \partial_{\theta_i}q= 0 \quad\hbox{ and  }\quad \sum_{1\le i<j\le n} \partial_{\theta_j}\partial_{\theta_i}q= 0 \Big\}~.
\end{equation}

Given $0\le k\le \lfloor \frac{n}{2}\rfloor$, a non-crossing pairing of length $k$ is a
list $(C_1, C_2, \ldots, C_k)$ with
\begin{align*}
 &C_r=(i_r,j_r) \text{ for $1\le i_r<j_r\le n$ for each $1 \leq r \leq k$ and,}\\
 &\qquad \text{either } i_r<j_r<i_s<j_s  \text{ or }  i_s<i_r<j_r<j_s\,\text{ for any $1\le r<s\le k$}.
\end{align*}
Given a non-crossing pairing $C=(C_1,C_2,\ldots C_k)$, we define
\begin{equation}\label{eq:Deltadef}
\Delta_C = \big(\theta_{j_1}-\theta_{i_1}\big)\big(\theta_{j_2}-\theta_{i_2}\big)\cdots \big(\theta_{j_k}-\theta_{i_k}\big)\,.
\end{equation}
Here $\Delta_C=1$ if $k=0$. Remark that $j_1<j_2<\cdots< j_k$.
The following proposition shows that there is a relationship between the non-crossing
partition condition and the differential equations from Equation \eqref{eq:defEQH}.

\begin{prop}\label{prop:harmelem}
The set
$${\mathcal D}'_n =\big\{ \Delta_C: \text{ $C=(C_1,C_2,\ldots, C_k)$  non-crossing pairing and $0\le k\le \lfloor \frac{n}{2}\rfloor$}\big\}
$$
is contained in $EQH_n $.
\end{prop}

\begin{proof}
To show that  $\Delta_C$ is  contained in $EQH_n $, we fix $C$. We need to show that $\Delta_C$ satisfies the differential equation conditions in Equation~\eqref{eq:defEQH}.

For the first  defining equation of $EQH_n$, we have
\begin{align*}
 \sum_{1\le i\le n} \partial_{\theta_i}\Delta_C&=\sum_{1\le r\le k}  (\partial_{\theta_{i_r}}+\partial_{\theta_{j_r}})\Delta_C + \sum_{i\notin \bigcup_{r=1}^k C_r}  \partial_{\theta_i}\Delta_C\\
 &=\sum_{1\le r\le k}  (\partial_{\theta_{i_r}}+\partial_{\theta_{j_r}})\Delta_C =0\,.
  \end{align*}
%For the first equality we simply reorganize the terms. For the second equality, we remark that if $i\notin \bigcup C$ then the variable $\theta_i$ does not appear in $\Delta_C$, so
%$ \partial_{\theta_i}\Delta_C=0$.
%%MZ: I'm taking that sentence out because I think it is clearer without stating it.
For the last equality, fix $1\le r\le k$ and note that
$\Delta_C =(-1)^{r-1}(\theta_{j_r}-\theta_{i_r}) q$ for some  polynomial $q$
and so for each $r$,
$$(\partial_{\theta_{i_r}}+\partial_{\theta_{j_r}})\Delta_C=(-1)^{r-1}
(\partial_{\theta_{i_r}}+\partial_{\theta_{j_r}})(\theta_{j_r}-\theta_{i_r})q =0~.
$$

For the second defining equation of $EQH_n $, we decompose the sum over pairs $1\le i<j\le n$ according to
whether (a) $|\{i,j\}\cap \bigcup_{r=1}^k C_r|<2$, (b) $C_r = (i,j)$ for some $r$, or
(c) $i,j$ appears in two different $C_{r},C_{s}$.

In case (a), if  $|\{i,j\}\cap \bigcup C|<2$, then one of $\theta_i$ or $\theta_j$ does not appear in $\Delta_C$ and we have $ \partial_{\theta_j} \partial_{\theta_i}\Delta_C=0$.

In case (b), we have that the product $\theta_{j_r}\theta_{i_r}$ does not appear in $\Delta_C$ and we also have $ \partial_{\theta_{j_r}} \partial_{\theta_{i_r}}\Delta_C=0$.

Thus we know that only case (c) contributes to the sum and we can thus write
\begin{align*}
 \sum_{1\le i<j\le n} \partial_{\theta_j}\partial_{\theta_i}\Delta_C
 &= \sum_{1\le r<s\le k} \sum_{i\in C_{r} \atop j\in C_{s}} \pm  \partial_{\theta_j}\partial_{\theta_i}\Delta_C~.
\end{align*}
In the second sum on the right hand side, we have to be careful as when we pick
$i\in C_{r}$ and $ j\in C_{s}$ we are not guaranteed that $i<j$ so a sign may be needed
in  order to keep the equality. We will make a careful study of all possibilities
for fixed $1\le r<s\le k$. First, we rearrange the terms of $\Delta_C$ to  bring the terms
$(\theta_{j_{r}}-\theta_{i_{r}})(\theta_{j_{s}}-\theta_{i_{s}})$ in front  performing $(r-1)+(s-2)$ anticommutations, we have
  $$\Delta_C  =(-1)^{r+s-1}(\theta_{j_{r}}-\theta_{i_{r}})(\theta_{j_{s}}-\theta_{i_{s}})q$$
  for some polynomial $q$. Remark that $i_{r},j_{r},i_{s},j_{s}$ satisfy either the inequalities
  $$  i_r<j_r<i_s<j_s \qquad \text{ or }\qquad  i_s<i_r<j_r<j_s.$$
 The only concern is their relative order, we can thus assume that we have the numbers $1,2,3,4$.
  There are two possibilities, $((i_r,j_r), (i_s,j_s))$ is equal to $((1,2),(3,4))$ or $((2,3),(1,4))$.
  In the first case we have
$$ \big( \partial_{\theta_3}\partial_{\theta_1} +  \partial_{\theta_3}\partial_{\theta_2} +  \partial_{\theta_4}\partial_{\theta_1} +  \partial_{\theta_4}\partial_{\theta_2}\big)
     (\theta_{2}-\theta_{1})(\theta_{4}-\theta_{3}) =0,
$$
 and in the second case we get
$$ \big( \partial_{\theta_2}\partial_{\theta_1} +  \partial_{\theta_4}\partial_{\theta_2} +  \partial_{\theta_3}\partial_{\theta_1} +  \partial_{\theta_4}\partial_{\theta_3}\big)
     (\theta_{3}-\theta_{2})(\theta_{4}-\theta_{1}) =0.
$$
And this shows that $\Delta_C\in EQH_n $ for all non-crossing pairings $C$.
\end{proof}

The set ${\mathcal D}'_n$ is not linearly independent, for example for $n=3$ and $k=1$, we have the following three non-crossing pairing:
$((1,2))$, $((1,3))$ and $((2,3))$, but
\[
\Delta_{((1,2))} - \Delta_{((1,3))} + \Delta_{((2,3))} =0 ~.
\]
We want to select a linearly independent subset of ${\mathcal D}'_n$. We proceed as follows:
consider a sequence
$\alpha = (a_1, a_2, \ldots, a_n) \in \{0, 1\}^n$
such that $\sum_{i=1}^r a_i \leq r/2$ for all $1 \leq r \leq n$.
Such sequences are known as \defn{ballot sequences}.
If ever it is the case that $\sum_{i=1}^r a_i > r/2$ then we say that
$\alpha$ \defn{breaks the ballot condition} at position $r$.

%In Section~\ref{sec:path} we will develop a more visual interpretation of
%ballot sequences and interpret them as paths that stay above the diagonal.
Given a ballot sequence $\alpha$ we build a non-crossing pairing $C(\alpha)$ by first replacing all $0$s
by open parentheses $0\mapsto$`(',
and all $1$s by close parentheses $1\mapsto$`)',
and then do the natural maximal pairing of parenthesis. The positions of the pairings
give us in lexicographic order a non-crossing pairing which we shall denote $C(\alpha)$.
Since $\alpha$ is a ballot sequence, every closed parenthesis is matched
and some open parentheses might remain unpaired.
The natural pairing of parenthesis guarantees that the result will be non-crossing. For example,
\[
\alpha=0010001101 \qquad\mapsto\qquad (()((())() \qquad\mapsto\qquad C(\alpha)=((2,3),(6,7),(5,8),(9,10)) \,.
\]

The total number of ballot sequences of size $n$ is equal to $\binom{n}{\lfloor{n/2}\rfloor}$
(see \cite[\href{https://oeis.org/A001405}{A001405}]{OEIS}).
The number of ballot sequences graded by the number of $1$'s in the sequence
(see \cite[\href{https://oeis.org/A008315}{A008315}]{OEIS})
is given in Figure \ref{table:ballotseq}.

\begin{figure}[h]
\begin{tabular}{c||cccccc}
$n=1$&1&&&&\\
$n=2$&1&1&&&\\
$n=3$&1&2&&&\\
$n=4$&1&3&2&&\\
$n=5$&1&4&5&&\\
$n=6$&1&5&9&5&\\
$n=7$&1&6&14&14&\\
$n=8$&1&7&20&28&14\\
$n=9$&1&8&27&48&42\\
\end{tabular}
\caption{The number of ballot sequences of length $n$ with exactly $k$ $1$s
with $1 \leq n \leq 9$ and $1 \leq k \leq \lfloor \frac{n}{2} \rfloor$.
These will be shown to be the graded dimensions of $EQH_n \simeq EQC_n$.}
\label{table:ballotseq}
\end{figure}

Given this construction we have the following Proposition.

\begin{prop}\label{prop:harmbasis}
The set
$${\mathcal D}_n =\big\{ \Delta_{C(\alpha)}:  \alpha \in \{0, 1\}^n \text{ a ballot sequence}\big\}
$$
is contained in $EQH_n$ and is linearly independent.
\end{prop}

\begin{proof}
  The first statement follows from Proposition~\ref{prop:harmelem} since ${\mathcal D}_n \subseteq {\mathcal D}'_n \subseteq EQH_n$.
To show the linear independence, fix $\alpha$ a ballot sequence and let $C(\alpha)=((i_1,j_1),\ldots,(i_k,j_k))$ be its non-crossing pairing.  We remark that the sequence of numbers
$j_1<j_2<\cdots<j_k$ corresponds to the position of the $1$s in $\alpha$.
Using the monomial ordering described in Section \ref{sec:linbasis}
and by inspection of the product in Equation~\eqref{eq:Deltadef},
we observe that the term $\theta_{j_1}\theta_{j_2}\cdots\theta_{j_k}$ is the smallest lexicographic
monomial in $\Delta_{C(\alpha)}$.
For different ballot sequences $\alpha$ we get different positions of the $1$s in
$\alpha$ and thus different smallest lexicographic monomials,
which shows the independence of ${\mathcal D}_n$.
\end{proof}

\begin{remark} For a fixed $0\le k\le \lfloor \frac{n}{2}\rfloor$, the set
$${\mathcal D}^{(k)}_n =\big\{ \Delta_{C(\alpha)}:  \alpha \in \{0, 1\}^n \text{ a ballot sequence with $k$  1s}\big\}$$
spans a subspace of $R_n$ of degree $k$.
It is known that the ballot sequences with $k$ 1s are in bijection with standard tableaux of shape $(n-k,k)$.
If  the variables $\theta$  were  commutative, the space spanned by ${\mathcal D}^{(k)}_n$ would be the same as the space spanned
by the Specht polynomials indexed by standard tableaux and therefore
would be an irreducible symmetric group module.
Here the situation appears to be related, but is in fact is quite different.

A small example is informative.
Consider $n=4$ and $k=2$.
There are two ballot sequences $0101$ and $0011$.
The associated two non-crossing pairings are
$((1,2),(3,4))$ and $((2,3),(1,4))$ and we have
$$\Delta_{C(0101)}=  (\theta_{2}-\theta_{1})(\theta_{4}-\theta_{3}) \quad\text{ and }\quad \Delta_{C(0011)}=(\theta_{3}-\theta_{2})(\theta_{4}-\theta_{1}).
$$
On the other hand, the two standard tableaux associated to $0101$ and $0011$ are
$$
T_1=\begin{tikzpicture}[scale=0.5,baseline=12pt]
	\draw[help lines] (0,0) grid (2,2);
	\node at (.5,.5){$\scriptstyle 1$};
	\node at (.5,1.5){$\scriptstyle 2$};
	\node at (1.5,.5){$\scriptstyle 3$};
	\node at (1.5,1.5){$\scriptstyle 4$};
\end{tikzpicture}
\qquad\text{ and }\qquad
T_2=\begin{tikzpicture}[scale=0.5,baseline=12pt]
	\draw[help lines] (0,0) grid (2,2);
	\node at (.5,.5){$\scriptstyle 1$};
	\node at (.5,1.5){$\scriptstyle 3$};
	\node at (1.5,.5){$\scriptstyle 2$};
	\node at (1.5,1.5){$\scriptstyle 4$};
\end{tikzpicture}\,.
$$
A standard construction of the symmetric group irreducible of shape  $(2,2)$ from the
tableaux $T_1$ and $T_2$ is to use the Garnir polynomials
\begin{align*}
  \Delta_{T_1} &=  (\theta_{2}-\theta_{1})(\theta_{4}-\theta_{3}) =\Delta_{C(0101)}\\
  \Delta_{T_2} &= (\theta_{3}-\theta_{1})(\theta_{4}-\theta_{2}).
\end{align*}
Unfortunately $ \Delta_{T_2}\notin EQH_n$.  In commutative variables, the  span of the $\{\Delta_{T_1},\Delta_{T_2}\}$ (an irreducible module) would be the same as
 the span of $\{\Delta_{C(0101)},\Delta_{C(0011)}\}$. But  for anticommutative variables, it is a  different story.\footnote{A more correct construction would be to apply the Young idempotent associated to $T_2$ to the monomial associated to $T_2$ using Hivert's action.
In this case  we get $\Delta'_{T_2}=  \theta_3\theta_4 -  \theta_1\theta_4 -  \theta_2\theta_3 +  \theta_1\theta_2 \notin EQH_n$. The span $\{\Delta_{T_1},\Delta'_{T_2}\}$ is a symmetric group irreducible module but is not fully contained in $EQH_n$.}
\end{remark}

%%%%%%%%%%%%%%%%%%%%%%%%%%%%%%%%%%%%%%
%%%%%%%%%%%%%%%%%%%%%%%%%%%%%%%%%%%%%%
%%%%%%%%%%%%%%%%%%%%%%%%%%%%%%%%%%%%%%
\section{A linear basis of the ring}\label{sec:linbasis}

Again let $n$ be a fixed positive integer and $R_n = {\mathbb Q}[\theta_1, \theta_2, \ldots, \theta_n]$.
We have thus far represented the basis for
$R_n$ as the elements $\theta_A$ with $A \subseteq [n]$.  Define $\alpha(A) \in \{ 0,1\}^n$ to be
the sequence $a_1 a_2 a_3 \cdots a_n$ with $a_i = 1$ if $i \in A$ and
$a_i = 0$ if $i \notin A$ so that
\[
\theta_A = \theta_1^{a_1} \theta_2^{a_2} \cdots \theta_n^{a_n} := \theta^{\alpha(A)}~.
\]
For such a sequence $\alpha \in \{0,1\}^n$, let $m_1(\alpha) := \sum_{i=1}^n a_i$
represent the number of $1$s in the string.  This will also be the degree of the monomial
$\theta^{\alpha}$.

For sequences $\alpha \in \{ 0, 1 \}^n$, define elements $G_\alpha$ by
\begin{equation}\label{eq:Gdef1}
G_{1^s0^{n-s}} = F_{1^s}
\end{equation}
and if $\alpha \neq 1^s 0^{n-s}$, then $\alpha$ is of the form $u01^s0^{n-k-s}$ for some string $u$
of length $k-1$ and we recursively define
\begin{equation}\label{eq:Gdef2}
G_{u01^s0^{n-k-s}} = G_{u1^s0^{n-k-s+1}} - (-1)^{m_1(u)} \theta_k G_{u1^{s-1}0^{n-k-s+2}}~.
\end{equation}

We will show below that the recurrence for the
$G_\alpha$ is defined so that they are $S$-polynomials \cite{CLO}
for elements of the ideal $I_n$.
In commutative variables, similar polynomials were defined by Aval-Bergeron-Bergeron~\cite{AB,ABB}
as a (complete) subset of $S$-polynomials needed to compute all possible
$S$-polynomials in the Buchburger algorithm for a Gr\"obner basis.
It is not given that one can easily describe such a set of $S$-polynomials and here we have adapted
the definition for working in the exterior algebra.

\begin{example} For $\alpha = 010110$ and $\beta = 001100$, we compute
the elements $G_\alpha$ and $G_\beta$ using the definition.
\begin{align*}
G_{010110} &= G_{011100} + \theta_3 G_{011000} = (G_{111000} - \theta_1 G_{110000})
+ \theta_3(G_{110000} - \theta_1 G_{100000})\\
&= \theta_2 \theta_4 \theta_5 + \theta_2 \theta_4 \theta_6 + \theta_2 \theta_5 \theta_6
+ 2 \theta_3 \theta_4 \theta_5 + 2 \theta_3 \theta_4 \theta_6 + 2 \theta_3 \theta_5 \theta_6
+ \theta_4 \theta_5 \theta_6
\end{align*}

and we have that
\begin{align*}
G_{001100} &= G_{011000} - \theta_2 G_{010000} = (G_{110000} - \theta_1 G_{100000})
- \theta_2 (G_{100000} - \theta_1 G_{000000})\\
&= \theta_3 \theta_4 + \theta_3 \theta_5 + \theta_3 \theta_6 + \theta_4 \theta_5 + \theta_4 \theta_6 + \theta_5 \theta_6.
\end{align*}
\end{example}

We follow \cite{CLO} for the convention of lexicographical ordering on monomials.
Given vectors $u,v$ with non-negative integer entries, we say $u < v$ lexicographically if
there exists an index $j \ge 1$ such that $u_{i} = v_{i}$ for all $1 \le i < j$ but
$u_{j} < v_{j}$.  Monomials of $R_{n}$ are ordered by their exponent vectors.  More
precisely, $\theta_{A} < \theta_{B}$ if $\alpha(A) < \alpha(B)$ lexicographically.  For example, we have
$\theta_{1} > \cdots > \theta_{n}$ and the lexicographically largest monomial in the above example
$G_{001100}$ is $\theta_{3}\theta_{4}$.  The latter demonstrates an important property of these
elements stated in the following proposition.
\begin{prop}\label{prop:largest}
The largest lexicographic term in $G_\alpha$ is $\theta^\alpha$.
\end{prop}

The proof of Proposition \ref{prop:largest} follows by induction
on the length of $\alpha$ and
from a lemma that is analogous to Lemma 3.3 of \cite{AB}.
The recursion in this result is really the origin of the definition
of $G_\alpha$ because Equation \eqref{eq:Gdef2}
was adapted so that this lemma holds.  It follows that
the set $\{ G_\alpha \}_{\alpha \in \{0,1\}^n}$ is a basis
for $R_n$.

The argument for the proposition is elementary (chasing the largest lexicographic
term in \eqref{eq:0start} and \eqref{eq:1start})
and so we do not include it, however the proof of the
following result comes from careful analysis of the
terms arising in the recursive definition of the $G_\alpha$.

\begin{lemma}\label{lem:LT}
 Let $\alpha \in \{0,1\}^{n-1}$, then
\begin{align}
G_{0\alpha} &= G_\alpha(\theta_2, \theta_3, \ldots, \theta_n) \hbox{ and }\label{eq:0start}\\
G_{1\alpha} &= \theta_1 G_{0\alpha} + P_{\alpha}(\theta_2, \theta_3, \ldots, \theta_n) \label{eq:1start}
\end{align}
for some polynomial
$P_{\alpha}(\theta_2, \theta_3, \ldots, \theta_n)
\in {\mathbb Q}[\theta_2, \theta_3, \ldots, \theta_n]$.
\end{lemma}

\begin{remark}\label{rem:shift}
By convention, the length of the index for our polynomials indicates in which polynomial space we are.
For example if $\beta   \in \{0,1\}^{n}$ then $G_\beta\in R_n$. For $\alpha \in \{0,1\}^{n-1}$  in Lemma~\ref{lem:LT}, when we write
$G_\alpha(\theta_2, \theta_3, \ldots, \theta_n)$ we mean $G_\alpha\in  R_{n-1}$
embedded in $R_n$ with the substitution
$\theta_i:=\theta_{i+1}$. A similar convention will be followed for $P_\alpha$.
\end{remark}

\begin{proof}[Proof of Lemma~\ref{lem:LT}]
The proof will proceed by induction on $n-i$ where $i$ is the number of trailing $0$s in $\alpha$.
The base case $0=n-n$ with $n$ zeros, is $0\alpha = 0^n$ and we have
\[
G_{0^n} = F_{1^0}(\theta_1,\theta_2 ,\ldots,\theta_n) = 1
= G_{0^{n-1}}(\theta_2, \theta_3, \ldots, \theta_n)~.
\]
We then consider the case  $0\alpha=01^s0^{n-s-1}$. The polynomials $F_{1^s}$ satisfy the  following identity
\begin{equation}\label{eq:relF}
F_{1^s}(\theta_1,\theta_2 ,\ldots,\theta_n)
= \theta_1 F_{1^{s-1}}(\theta_2 ,\ldots,\theta_n) + F_{1^s}(\theta_2,\theta_3 ,\ldots,\theta_n).
\end{equation}
This follows directly from the definition~\eqref{eq:defF} where we split the sum in two parts depending if $1\in A$ or not.
The definition of $G_{01^s0^{n-s-1}}$  gives us
\begin{align*}
	G_{01^s0^{n-s-1}}&= G_{1^s0^{n-s}} - \theta_1 G_{1^{s-1}0^{n-s+1}}
	  =F_{1^s}(\theta_1,\theta_2 ,\ldots,\theta_n) - \theta_1 F_{1^{s-1}}(\theta_2 ,\ldots,\theta_n) \\
	  &= F_{1^s}(\theta_2,\theta_3 ,\ldots,\theta_n)= G_{\alpha}(\theta_2,\theta_3 ,\ldots,\theta_n)\,.
\end{align*}

To finish the  proof of Equation~\eqref{eq:0start} by induction,
let us assume that $\alpha$ is not of the form
$01^s0^{n-s-1}$ for some $s>0$.
Instead we have $0\alpha=0w01^s0^{n-k-s}$ for some $s>0$ and some string $w$
of length $k-2$.
For $0\alpha=0w01^s0^{n-k-s}$, we have $n-k-s$ trailing zeros.
Remark that for $0w1^s0^{n-k-s+1}$ and $0w1^{s-1}0^{n-k-s+2}$
we have more trailing zeros than that of $0\alpha$ and we will use the
induction hypothesis with~\eqref{eq:0start} in the
equality~\eqref{eq:IH0} below.%%%%%%%
\begin{align}
	G&_{0w01^s0^{n-k-s}}= G_{0w1^s0^{n-k-s+1}} -  (-1)^{m_1(0w)} \theta_k G_{0w1^{s-1}0^{n-k-s+2}} \nonumber \\
	 & = G_{w1^s0^{n-k-s+1}}(\theta_2,\theta_3 ,\ldots,\theta_n) - (-1)^{m_1(0w)} \theta_k G_{w1^{s-1}0^{n-k-s+2}}(\theta_2,\theta_3 ,\ldots,\theta_n) \label{eq:IH0}\\
	 & = \big[G_{w1^s0^{n-k-s+1}} - (-1)^{m_1(w)} \theta_{k-1} G_{w1^{s-1}0^{n-k-s+2}}\big](\theta_2,\theta_3 ,\ldots,\theta_n)\label{eq:detail0}\\
	 & = G_{w01^s0^{n-k-s}} (\theta_2,\theta_3 ,\ldots,\theta_n)
	   = G_{\alpha} (\theta_2,\theta_3 ,\ldots,\theta_n).  \label{eq:final0}
\end{align}
In~\eqref{eq:detail0} the expression inside the square bracket $[\cdots]$ is treated as a polynomial in the variables $\theta_1,\ldots,\theta_{n-1}$ in $R_{n-1}$ (see Remark~\ref{rem:shift}).
Hence, the variable $\theta_{k}$ from \eqref{eq:IH0} must
be replaced by $\theta_{k-1}$ in \eqref{eq:detail0}.
Also, $m_1(0w)=m_1(w)$. The expression we get is exactly the definition of
$G_{w01^s0^{n-k-s}} \in R_{n-1}$ and Equation~\eqref{eq:final0} follows.
This concludes the proof of \eqref{eq:0start}.

We next prove Equation \eqref{eq:1start} by induction.  The base case is if $1\alpha = 1^{s+1}0^{n-s-1}$,
then using Equation \eqref{eq:relF} we have
\begin{align}
G_{1^{s+1}0^{n-s-1}} &= F_{1^{s+1}}(\theta_1,\theta_2 ,\ldots,\theta_n) \nonumber\\
&= \theta_1 F_{1^s} (\theta_2, \theta_3, \ldots, \theta_n)+ F_{1^{s+1}}(\theta_2, \theta_3, \ldots, \theta_n)\nonumber \\
&= \theta_1 G_{01^s0^{n-s-1}} + F_{1^{s+1}}(\theta_2, \theta_3, \ldots, \theta_n)~.\label{eq:base1}
\end{align}
In~\eqref{eq:base1}, we use \eqref{eq:0start} with  $G_{01^s0^{n-s-1}}=G_{1^s0^{n-s-1}} (\theta_2,  \ldots, \theta_n)=F_{1^s} (\theta_2, \ldots, \theta_n)$.
We then let
$P_{1^{s}0^{n-s-1}} = F_{1^{s+1}}(\theta_1,\ldots,\theta_{n-1})$ and this shows that \eqref{eq:1start} holds in this case.

We now assume that $1\alpha\ne 1^{s+1}0^{n-s-1}$. Therefore $1\alpha=1w01^s0^{n-k-s}$ for some string $w$
of length $k-2$. We have
\begin{align}
	G&_{1w01^s0^{n-k-s}}= G_{1w1^s0^{n-k-s+1}} -  (-1)^{m_1(1w)} \theta_k G_{1w1^{s-1}0^{n-k-s+2}} \nonumber \\
	 & = \big(\theta_1G_{0w1^s0^{n-k-s+1}} +P_{w1^s0^{n-k-s+1}}(\theta_2,\theta_3 ,\ldots,\theta_n)\big) \label{eq:IH1}\\
	 &\qquad\qquad  - (-1)^{m_1(1w)} \theta_k  \big(\theta_1G_{0w1^{s-1}0^{n-k-s+2}} +P_{w1^{s-1}0^{n-k-s+2}}(\theta_2,\theta_3 ,\ldots,\theta_n)\big) \nonumber\\
	 & = \theta_1\big(G_{0w1^s0^{n-k-s+1}}   - (-1)^{m_1(w)}  \theta_k G_{0w1^{s-1}0^{n-k-s+2}}\big)	 \label{eq:detail1}\\
	 &\qquad\qquad  + \big[P_{w1^s0^{n-k-s+1}}  - (-1)^{m_1(1w)} \theta_{k-1} P_{w1^{s-1}0^{n-k-s+2}}\big] (\theta_2,\theta_3 ,\ldots,\theta_n). \nonumber
%	 & = \theta_1 G_{0\alpha} + P_\alpha(\theta_2, ,\theta_3 ,\ldots,\theta_n)  \nonumber
\end{align}
In \eqref{eq:IH1}, we have used the induction hypothesis of~\eqref{eq:1start} on both terms.
In  \eqref{eq:detail1}, we group together the terms with $\theta_1$ in front, using the identity
$(-1)^{m_1(1w)} \theta_k \theta_1 = (-1)^{m_1(1w)+1} \theta_1 \theta_k  = (-1)^{m_1(w)} \theta_1 \theta_k $.
The term with $\theta_1$ in \eqref{eq:detail1} is the definition of $G_{0\alpha}$.
The expression inside the square bracket is a polynomial in $R_{n-1}$ that we take as the definition for $P_\alpha$.
This shows by induction that \eqref{eq:1start}
holds in all cases and concludes the  proof of the lemma.
\end{proof}

%%%%%%%%%%%%%%%%%%%%%%%%%%%%%%%%%%%%%%
%%%%%%%%%%%%%%%%%%%%%%%%%%%%%%%%%%%%%%
%%%%%%%%%%%%%%%%%%%%%%%%%%%%%%%%%%%%%%
\section{A basis for the quotient}\label{sec:ballotbasis}

%To each sequence $\alpha \in \{0, 1\}^n$, we associate a path starting at the origin
%and extending into the first quadrant of the $x,y$-plane.  The $i^{th}$ step of this
%path will be a unit in the $(1,0)$-direction if $a_i =1$ and it will be a unit in the $(0,1)$-direction
%if $a_i = 0$.  We say that the sequence $\alpha$ \defn{crosses the diagonal} if
%there is a point on the path which lies at $(a,a)$ and the next step is in the $(1,0)$
%direction.  Otherwise we say that $\alpha$ \defn{stays above the diagonal}.  Note that
%$\alpha$ stays above the diagonal if $\sum_{i=1}^r a_i \leq r/2$
%for all $1 \leq r \leq n$
%and it crosses the diagonal otherwise. The ballot sequences we encountered earlier are the one that  stay above the diagonal.
%\mike{I think that we should state everything in terms of ballot sequences
%or paths, but not both.  I see no obvious reason to switch between the two models
%and we should make the presentation as simple as possible.}
%
%\begin{example} Let $n=6$ and $\alpha = 010110$ and $\beta = 001100$, then the corresponding
%paths are
%\begin{center}
%\begin{tikzpicture}[scale=.75]
%  \draw[dotted] (0, 0) grid (4, 4);
%  \draw[rounded corners=1, color=black, line width=2] (0, 0) -- (0, 1) -- (1, 1) -- (1, 2) -- (2, 2) -- (3, 2) -- (3, 3);
%\end{tikzpicture}
%\hskip .5in
%\begin{tikzpicture}[scale=.75]
%  \draw[dotted] (0, 0) grid (4, 4);
%  \draw[rounded corners=1, color=black, line width=2] (0, 0) -- (0, 1) -- (0, 2) -- (1, 2) -- (2, 2) -- (2, 3) -- (2, 4);
%\end{tikzpicture}
%\end{center}
%\end{example}

The elements $G_\alpha$ are defined so that we could use them to
identify a nice basis of the ideal $I_n$.  Our first result establishes
that the $G_\alpha$ such that $\alpha$ is not a ballot sequence
are in the ideal.  The slightly more difficult step is to show
that these elements also span the ideal.

\begin{prop}\label{prop:notballotimpliescontains}
If $\alpha \in \{0,1\}^n$ is not a ballot sequence, then $G_\alpha \in I_n$.
\end{prop}

\begin{proof}
A sequence $\alpha \in \{ 0, 1\}^n$ is either of the form
$\alpha = 1^s0^{n-s}$ for some $s > 0$ or
$\alpha = u01^s0^{n-s-k}$ for some $s>0$ and some $u \in \{0,1\}^{k-1}$.

In the first case, $\alpha$ breaks the ballot condition in position $1$
and by Equation \eqref{eq:Gdef1}, $G_{1^s0^{n-s}} = F_{1^s}$
is in the ideal $I_n$.

Now the other case is established by
induction on the position of the last $1$ in $\alpha$.
We assume that
$\alpha = u01^s0^{n-s-k}$ and, by Equation \eqref{eq:Gdef2},
$G_\alpha$ is in $I_n$ if both
$G_{u1^s0^{n-k-s+1}}$ and $G_{u1^{s-1}0^{n-k-s+2}}$
are elements of $I_n$.

Assume that $u01^s0^{n-s-k}$
breaks the ballot condition for the first time at position $r$.
If $r<k$, then
$u1^s0^{n-k-s+1}$ and $u1^{s-1}0^{n-k-s+2}$
both break the ballot condition also at position $r$.
Since $\alpha_k=0$, $\alpha$ does not break the ballot condition for the first
time at $r=k$, so the other possibility is that is $r>k$.
In this case $u01^{r-k}$ with $r-k \leq s$ breaks the ballot condition for the first time
and therefore so does $u1^{r-k-1}$ and so do both
$u1^s0^{n-k-s+1}$ and $u1^{s-1}0^{n-k-s+2}$.
By our inductive hypothesis this implies $G_\alpha \in I_n$.

Therefore by induction, $\alpha$ breaks the ballot condition implies $G_\alpha \in I_n$
for all $\alpha \in \{0,1\}^n$.
\end{proof}

We will  show that the ideal lies in the span of the $G_\alpha$ such
that $\alpha$ breaks the ballot condition therefore establishing our main theorem.

\begin{theorem}\label{thm:basisofideal}
The set $A_n:=\big\{ G_\alpha : \alpha \in \{0,1\}^n \text{ \it breaks the ballot condition}\big\}$
is a $\mathbb Q$--linear basis of the ideal $I_n$.
\end{theorem}

The proof of this theorem uses our understanding of the harmonic space $EQH_n\cong EQC_n$.
In Proposition~\ref{prop:harmbasis} we found that $\dim(EQH_n)=\dim(EQC_n)$ is at least the number of ballot sequences.
We first establish a small lemma about a spanning set for the quotient $EQC_n$ showing that  the dimension is at most  the number of ballot sequences.
Therefore we have equality
and the set ${\mathcal D}_n$ in Proposition~\ref{prop:harmbasis} is in fact a basis
of $EQH_n$.

\begin{lemma}\label{lem:spanofquotient}
The set $B_n=\big\{ \theta^\beta : \beta \in \{0,1\}^n \text{ \it is a ballot sequence}\big\}$
$\mathbb Q$--spans the quotient $R_n\big/I_n$.
\end{lemma}

\begin{proof}
Order the monomials lexicographically and let $\theta^\gamma$ be the smallest monomial
that is not in the $\mathbb Q$--span of $B_n$
(modulo  $I_n$). We must have  that $\gamma$ breaks the ballot condition,
since otherwise $\theta^\gamma\in  B_n$. Therefore, Proposition~\ref{prop:notballotimpliescontains}
tells us  that $G_\gamma \in I_n$. Proposition~\ref{prop:largest} says that $G_\gamma = \theta^\gamma + \sum_{\beta<\gamma}  c_\beta \theta^\beta$. Hence, modulo $I_n$, we have
$$ \theta^\gamma \equiv \theta^\gamma - G_\gamma = - \sum_{\beta<\gamma} c_\beta \theta^\beta\,.$$
The right hand side is a linear combination of monomials strictly smaller than $\theta^\gamma$.
By the choice of $\theta^\gamma$, all such monomials are in the $\mathbb Q$--span of $B_n$.
 Therefore $\theta^\gamma$ is also in the $\mathbb Q$--span of $B_n$, a contradiction. We must conclude there are no such $\theta^\gamma$ and all monomials are in the $\mathbb Q$--span of $B_n$ modulo the ideal $I_n$.
\end{proof}

\begin{proof}[Proof of Theorem~\ref{thm:basisofideal}]
Let $d_n$ be the number of ballot sequences of size $n$.
We have
$$ \dim EQH_{n} \le d_n \le \dim EQC_{n},$$
where the first inequality follows from Proposition~\ref{prop:harmbasis} and second follows from Lemma~\ref{lem:spanofquotient}. By Proposition~\ref{prop:EQC_EQH}, we have $d_{n} = \dim EQC_{n}$.
Let ${\mathbb Q}A_n$ be the ${\mathbb Q}\text{-span}$ of the elements of $A_n$.
Similarly, let ${\mathbb Q}B'_n$ be the ${\mathbb Q}\text{-span}$ of the set $\big\{ G_\beta : \beta \in \{0,1\}^n \text{ \it is a ballot sequence}\big\}$.
Using Proposition~\ref{prop:largest} we have that
$$R_n ={\mathbb Q}B'_n \oplus  {\mathbb Q}A_n
$$
Since ${\mathbb Q}A_n \subseteq I_n$ and $\dim {\mathbb Q}B'_n = d_n$, we conclude that ${\mathbb Q}A_n = I_n\,.$
\end{proof}

There are several straightforward consequence of this theorem which we state here.

\begin{cor}\label{cor:hilb}  The number of ballot sequences of length $n$ with $k$ entries $1$ is known~\cite{Stanley} to be equal to the number $f^{(n-k,k)}$ of standard tableaux of shape $(n-k,k)$.
Therefore we have  the Hilbert series of $EQH_n$ is
$$  \text{Hilb}_{EQH_n}(q) =  \sum_{k=0}^{\lfloor{n/2}\rfloor} f^{(n-k,k)} q^k\,.$$
\end{cor}

\begin{cor} The set ${\mathcal D}_n$ is a basis of $EQH_n$ and of $EQC_n$.
\end{cor}

\begin{cor}\label{cor:minimalGB}
The set $A_n$ is a (non-reduced, non-minimal) Gr\"obner basis of $I_n$. A minimal Gr\"obner basis for $I_n$ is given by
$$ \big\{ G_\alpha : \alpha \in \{0,1\}^n \text{ \it breaks the ballot condition only
at the rightmost  1 of }  \alpha\big\}\,.
$$
\end{cor}

\begin{remark} In this paper we have adopted the language of ballot sequences. An alternative (as in~\cite{B,AB}) is to use north-east lattice paths in the first quadrant from $(0,0)$ with
a north step for every $0$ and an east step  for every  $1$ as we read in a $0-1$ sequence. In such representation, a sequence is ballot if and only if it remains above
the diagonal.
\end{remark}

\end{document}